# ASYMPTOTIC THEORY FOR THE COX MODEL WITH MISSING TIME-DEPENDENT COVARIATE

BY JEAN-FRANÇOIS DUPUY, ION GRAMA AND MOUNIR MESBAH

*Université Toulouse 3, Université de Bretagne-sud and Université Paris 6*

The relationship between a time-dependent covariate and survival times is usually evaluated via the Cox model. Time-dependent covariates are generally available as longitudinal data collected regularly during the course of the study. A frequent problem, however, is the occurence of missing covariate data. A recent approach to estimation in the Cox model in this case jointly models survival and the longitudinal covariate. However, theoretical justification of this approach is still lacking. In this paper we prove existence and consistency of the maximum likelihood estimators in a joint model. The asymptotic distribution of the estimators is given along with a consistent estimator of the asymptotic variance.

**1. Introduction.** The commonly used Cox [6] regression model postulates that the hazard function for the failure time $T$ associated with a time-varying covariate $Z$ takes the form

$$\lambda(t; Z) = \lambda_0(t) \exp[\beta_0 Z(t)], \tag{1.1}$$

where $\beta_0$ is an unknown regression parameter and $\lambda_0$ is an unspecified baseline hazard function. The statistical problem is that of estimating $\beta_0$ and the cumulative baseline hazard function $\Lambda_0(t) = \int_0^t \lambda_0(s)\,ds$ on the basis of $n$ possibly right-censored survival times $X_1, \ldots, X_n$ and the corresponding covariates $Z_1, \ldots, Z_n$, where $Z_i$ is observed on the interval $[0, X_i]$.

By maximizing the partial likelihood [7], one can obtain an estimator of $\beta_0$ that is consistent and asymptotically normal with a covariance matrix which can be consistently estimated [1]. Letting $\hat{\beta}$ be the maximum partial likelihood estimator of $\beta$, Breslow [3, 4] suggested estimating $\Lambda_0(t)$ by

$$\hat{\Lambda}(t) = \sum_{X_i \leq t} \frac{\Delta_i}{\sum_{j=1}^{n} \exp[\hat{\beta} Z_j(X_i)] \mathbb{1}_{\{X_i \leq X_j\}}}.$$

Andersen and Gill [1] have shown weak convergence of the process $n^{1/2}(\hat{\Lambda} - \Lambda_0)$ to a Gaussian process.

To apply this methodology, one needs the knowledge of $\{Z(s) : 0 \leq s \leq t\}$ for all values $t \leq X$. This is generally not available. Common problems in survival analysis are presence of covariate measurement error (see among others Dafni and Tsiatis [8], Dupuy [12], Li and Lin [19], Tsiatis and Davidian [29], DeGruttola, Tsiatis and Wulfsohn [31] and Wulfsohn and Tsiatis [34]) and occurrence of missing covariate data (see [13, 14, 15, 20, 26, 35]).

A recent approach to estimation in the Cox model with a missing or mismeasured covariate consists in jointly modeling survival and the longitudinal covariate data. An extensive literature has now contributed to the estimation in such models (see [30] for a review and numerous references). However, rigorous proofs of the large-sample properties of estimators obtained from joint models remain an open problem. Note that simulations by Tsiatis and Davidian [29] show that the joint modeling approach should yield a consistent and asymptotically normal estimator for the regression parameter $\beta_0$. Li and Lin [19] provide simulations which also seem to point to the asymptotic validity of this approach for estimating parameters in the frailty model with covariate measurement error.

In this paper we propose a joint model for estimating parameters in the Cox model with missing values of a longitudinal covariate. Estimation in this joint model is carried out via nonparametric maximum likelihood (NPML) estimation. We prove consistency and asymptotic normality of the NPML estimator, and we give a consistent estimator for the limiting variance.

The paper is organized as follows: in Section 2 we describe the joint model and derive the likelihood function. In Section 3 we investigate the theoretical properties of the model, including identifiability and the existence of the NPML estimator. In Section 4 we show that the NPML estimator is consistent and asymptotically normal and we give a consistent estimator of its asymptotic variance.

**2. The statistical model and construction of the joint likelihood.** Suppose that $n$ subjects are observed. For each individual, we observe survival and covariate data. Denote by $T_i$ the random survival time for individual $i$.

We assume that survival is subject to right censoring, that is, instead of $T_i$, we actually observe $X_i = \min(T_i, C_i)$ and a failure indicator $\Delta_i = \mathbb{1}_{\{T_i \leq C_i\}}$, where $C_i$ is a random censoring time.

We examine the case of a single covariate $Z$ that is measured over time at the instants $0 = t_0 < t_1 < t_2 < \cdots$. We denote $Z_i(t_j)$ by $Z_{i,j}$. For $t > 0$, let $a_t = \max(k : t_k < t)$ be the index of the last observed value of $Z$ before time $t$.

The problem is as follows. Suppose that the data consist of i.i.d. replicates $(X_i, \Delta_i, Z_i(\cdot))$ $(i = 1, \ldots, n)$ of $(X, \Delta, Z(\cdot))$. For each subject $i$, the actually observed data is an incomplete random vector $\mathbf{Y}_i = (X_i, \Delta_i, Z_{i,0}, \ldots, Z_{i,a_{X_i}})$, where the covariate value at the time of failure $Z_i(X_i)$ is missing. The goal is to estimate the unknown true regression parameter $\beta_0$ and the cumulative hazard function $\Lambda_0(t) = \int_0^t \lambda_0(u) \, du$ $(t \geq 0)$ using the incomplete vectors $\mathbf{Y}_i$ $(i = 1, \ldots, n)$.

This work was motivated by a study whose design called for repeated measurements of a covariate to be made at different times $t_0 < t_1 < t_2 < \cdots$ on patients until drop-out. The objective is to evaluate the relationship between drop-out and the longitudinal values. The covariate being measured at the prespecified times $t_j$, Cox regression using model (1.1) is complicated by missingness of the covariate values at drop-out times. Some recent approaches to this problem [20, 26, 35] consist in extrapolating $Z_i(u)$ at failure time using available longitudinal data. However, for these methods to be valid, it is assumed that drop-out is nonignorable, that is, the probability of drop-out does not depend on the unobserved covariate value. This hypothesis does not hold in our setting, since hazard of drop-out at time $t$ depends on the unobserved $Z(t)$. We then propose to jointly model survival and the covariate in order to use full data available to estimate the parameters. Applications of this approach in psychometry and AIDS clinical trials can be found in [14, 15, 22], along with comparisons with alternative methods.

In order to derive asymptotic results for our estimators, we assume throughout that the following conditions C1–C7 are satisfied:

C1. Let $\tau$ be a finite time point at which any individual still under study is censored. Assume that $P(X \geq \tau) > 0$.
C2. Conditional on the observed path of the longitudinal covariate, the hazard function for $T_i$ is given by $\lambda_0(t) \exp[\beta_0 Z(t)]$.
C3. The covariates $Z_{i,j}$ have uniformly bounded total variation, namely, $\int_0^\infty |dZ_{i,j}(t)| + |Z_{i,j}(0)| \leq c$ for some finite $c > 0$ and all $i, j$.
C4. Let $f$ denote the joint density of $(Z_0, \ldots, Z_{a_t}, Z(t))$. Suppose that $f$ depends on an unknown parameter $\alpha$ $(\alpha \in \mathbb{R}^p)$, that $f$ is continuous with respect to $\alpha$ and has continuous second-order derivatives with respect to $\alpha$. Suppose also that $f$ is bounded and that, for any $t$, $f(z_0, \ldots, z_{a_t}, z(t); \alpha) = f(z_0, \ldots, z_{a_t}, z(t); \alpha')$ a.e. implies $\alpha = \alpha'$.
C5. The parameters $\alpha$ and $\beta$ are interior points of known compact sets $A \subset \mathbb{R}^p$ and $B \subset \mathbb{R}$, respectively. $\Lambda$ belongs to the set $L$ of absolutely continuous [with respect to the Lebesgue measure on $[0, \infty)$], nondecreasing functions $\Lambda$ such that $\Lambda(0) = 0$. Assume $\Lambda(\tau) < \infty$.

C6. Let $\theta = (\alpha, \beta, \Lambda)$, and note by $\theta_0 = (\alpha_0, \beta_0, \Lambda_0)$ the true value of $\theta$. Let $\Theta$ denote the parameter space $A \times B \times L$, and suppose that $\theta_0 \in \Theta$. Denote by $E_{\theta_0}[\cdot]$ the expectation of random variables taken under the true parameter. Suppose that $E_{\theta_0}[e^{\beta_0 Z(u)} \mathbb{1}_{\{u \leq X\}}]$ is bounded away from 0 on $[0, \tau]$, that $E_{\theta_0}[\int_0^X \{Z(u)\}^2 e^{\beta_0 Z(u)} \, d\Lambda_0(u)] > 0$, and that $-E_{\theta_0}[\frac{\partial^2}{\partial \alpha \, \partial \alpha^T} \ln f(Z_0, \ldots, Z; \alpha_0)]$ is positive definite.

C7. It is assumed that $T$ and $C$ are independent given the covariate $Z$. Moreover, we assume that the censoring distribution does not depend on the unobserved covariate value, or on $\theta$.

Condition C1 is a standard assumption that supposes that some individuals are at risk at the end $\tau$ of the experiment.

Condition C2 assumes for ease of presentation that hazard of failure at time $t$ depends on the time-varying covariate through its value at $t$. This could be relaxed, for example, by including a value $Z(t-h)$ $(h>0)$ (such as in $\lambda_0(t) \exp[\gamma_0 Z(t-h) + \beta_0 Z(t)])$ to study whether the variation in $Z$ between $(t-h)$ and $t$ influences survival. We shall note, however, that, in this case, $\beta_0$ and $\gamma_0$ are not identifiable if $Z$ is a time-independent covariate. We refer to Chen and Little [5] for their work on Cox regression with a missing time-independent covariate.

Condition C4 allows several kinds of parametric models to be used for the time-dependent covariate. For example, for each individual $i$, the $z_{i,j}$'s may be treated as a realization of a multivariate normal random vector, whose mean may possibly depend on explanatory variables (times of measurements $t_j$, treatment arms, covariates measured at the entry such as age, gender, ...). Various correlation structures may be assumed to take account of the correlation between measurements within each individual (see [10], Chapters 4 and 5). The parameter $\alpha$ would separate here into components for the mean and covariance structures. One may in this case impose additional conditions to ensure identifiability of covariance parameters, such as a minimum number of repeated measurements on some subjects. One may also use transition models (see [10], Chapters 7 and 10), where the conditional distribution of each $Z_{i,j}$ is modeled as a function of past responses $Z_{i,j-1}, \ldots$ and explanatory variables. Dupuy and Mesbah [14, 15] propose a joint model which uses a transition model for the longitudinal data. We refer to [10] for a detailed exposition of various parametric models for longitudinal data.

Condition C6 will ensure invertibility of a Fisher information operator in the proof of asymptotic normality of the estimators in the joint model.

Condition C7 is the usual condition of independent and noninformative censoring. It is usually satisfied in applications, in particular, when a subject is censored at $\tau$.

The probability measure induced by the observed $\mathbf{Y}$ is denoted by $P_\theta(d\mathbf{y}) = f_\mathbf{Y}(\mathbf{y}; \theta) \, d\mathbf{y}$ $(\theta \in \Theta)$. We shall obtain the likelihood $f_\mathbf{Y}(\mathbf{y}; \theta)$ for the vector of observations $\mathbf{y} = (x, \delta, z_0, \ldots, z_{a_x})$ by first writing the density of $(\mathbf{y}, z)$ for some value $z$ of $Z(X)$, and then by integrating over $z$. Actually, only a partial likelihood $L(\theta)$ for $\mathbf{Y}$ is specified, by discarding from $f_{\mathbf{Y}}(\mathbf{y}; \theta)$ terms adhering to censoring. From assumption C7, this does not influence maximization. The resulting likelihood $L(\theta)$ is

$$(2.2) \quad \int \{\lambda(x)\}^\delta \exp\left[\delta\beta z - \int_0^x \lambda(u) e^{\beta z(u)}\, du\right] f(z_0, \ldots, z_{a_x}, z; \alpha)\, dz,$$

where integration is over $z$ and the indetermination $0^0$ is set to be equal to 1. In the following, we shall denote by $l(\mathbf{y}, z; \theta)$ the integrand of (2.2).

**3. Nonparametric maximum likelihood estimation.** We shall first demonstrate identifiability of the proposed joint model. The proof is given in the Appendix.

PROPOSITION 1. *Under conditions* C1–C7, *the model is identifiable, that is, $L(\theta) = L(\theta')$ for almost all $\mathbf{y}$ implies $\theta = \theta'$, for $\theta, \theta' \in \Theta$.*

The problem of estimating $\theta$ is semiparametric, since the component $\Lambda$ is a function. Note that the maximum in $\Lambda$ of this likelihood function does not exist, so the principle of maximum likelihood is not applicable here. Nevertheless, this principle can be conveniently modified to yield a reasonable estimator of the function $\Lambda$, as well as of $\beta$ and $\alpha$.

We assume that there are no tied event times and that the number of events $p(n)$ increases with the sample size $n$. We reorder the indices of the data such that $X_1 < \cdots < X_{p(n)}$ $[p(n) \leq n]$ represent the increasingly ordered event times and $X_{p(n)+1} \leq \cdots \leq X_n$ represent the nondecreasingly ordered censoring times. To define an estimator of $\Lambda$ out of the likelihood (2.2), we proceed by the method of sieves [16], which consists in replacing the parameter space $\Theta$ by an appropriate approximating space $\Theta_n$ called a sieve (we refer to Li and Lin [19], McKeague [21] and Murphy and Sen [25] among others, for use of sieves in various settings of survival analysis). Precisely, instead of the functions $\Lambda = \Lambda(t), t \geq 0$, one considers increasing stepwise versions $\Lambda_n = \Lambda_n(t), t \geq 0$, of them with the unknown deterministic values $\Lambda_n(X_i) = \Lambda_{n,i}$ in the points $X_i, i = 1, \ldots, p(n)$. The sieve $\Theta_n$ is then

$$\{\theta = (\alpha, \beta, \Lambda_n) : \alpha \in \mathbb{R}^p, \beta \in \mathbb{R}, \Lambda_{n,1} \leq \cdots \leq \Lambda_{n,p(n)}, \Lambda_{n,i} \in \mathbb{R}, i = 1, \ldots, p(n)\}.$$

We shall estimate the values $\Lambda_{n,i}$ $[i = 1, \ldots, p(n)]$ and the parameters $\beta$ and $\alpha$ by maximizing the likelihood (2.2) over the parameter space $\Theta_n$, which means maximizing the pseudo likelihood

$$(3.3) \quad L_n(\theta) = \prod_{i=1}^n L^{(i)}(\theta)$$

obtained by multiplying over uncensored subjects $i$ $[i = 1, \ldots, p(n)]$ the individual contribution $L^{(i)}(\theta)$,

$$\int \Delta\Lambda_{n,i} \exp\left[\beta z - \sum_{k=1}^{p(n)} \Delta\Lambda_{n,k} e^{\beta z_i(x_k)} \mathbb{1}_{\{x_k \leq x_i\}}\right] f(z_{i,0}, \ldots, z_{i,a_{x_i}}, z; \alpha) \, dz,$$

and over censored subjects $i$ $[i = p(n) + 1, \ldots, n]$ the individual contribution

$$L^{(i)}(\theta) = \int \exp\left[-\sum_{k=1}^{p(n)} \Delta\Lambda_{n,k} e^{\beta z_i(x_k)} \mathbb{1}_{\{x_k \leq x_i\}}\right] f(z_{i,0}, \ldots, z_{i,a_{x_i}}, z; \alpha) \, dz,$$

where $\Delta\Lambda_{n,k} = \Delta\Lambda_n(X_k) = \Lambda_{n,k} - \Lambda_{n,k-1}$ $[k = 2, \ldots, p(n)]$ and $\Delta\Lambda_{n,1} = \Lambda_{n,1}$. The resulting estimator is usually referred to as nonparametric maximum likelihood estimator (NPMLE). We will use this terminology in the following, keeping in mind that the only part which is really nonparametric is just the representation of the baseline hazard. We next demonstrate that the estimator in the joint model does share the useful properties of the estimators from parametric models. We refer to [19, 23, 24, 28, 34] for use on NPMLE in various situations.

Proposition 2 shows that such an estimator exists in the proposed joint model. The proof is given in the Appendix.

PROPOSITION 2. *A maximizer $\hat{\theta}_n = (\hat{\alpha}_n, \hat{\beta}_n, \hat{\Lambda}_n)$ of $L_n(\theta)$ over $\theta \in \Theta_n$ exists and is achieved.*

To maximize the logarithm of the likelihood $L_n(\theta)$, we use the well-known approach used for the so-called expectation–maximization (EM) algorithm [9]. The rationale for this approach is that direct maximization of the integrated likelihood (3.3) is difficult, and that in the present setting, the maximizer $\hat{\theta}_n$ can be more easily characterized from an alternative EM-loglikelihood.

The following proposition provides an important characterization of the maximizer of $\sum_{i=1}^n \ln f_\mathbf{Y}(\mathbf{y}_i; \theta)$ [or, equivalently, of $L_n(\theta)$] on $\Theta_n$. This characterization will serve for the proof of asymptotic properties.

In the following, for any random variable $X$ with density function $f_X(x; \theta)$, we shall denote by $E_\theta[g(X)]$ the expected value of $g(X)$. Moreover, if $X$ and $Y$ are random variables, we shall denote by $E_\theta[g(X)|y]$ the expectation of $g(X)$ taken with respect to the conditional density function $f_{X|Y}(x|y; \theta)$ of $X$ given $Y = y$.

PROPOSITION 3. *The NPML estimator $\hat{\theta}_n$ satisfies the equation*

$$\hat{\Lambda}_n(t) = \int_0^t \frac{dH_n(u)}{W_n(u; \hat{\theta}_n)},$$

*where*

$$H_n(u) = n^{-1}\sum_{i=1}^{n}\Delta_i\mathbb{1}_{\{X_i\leq u\}} \quad \text{and} \quad W_n(u;\theta) = n^{-1}\sum_{i=1}^{n}E_\theta[e^{\beta Z(u)}\mathbb{1}_{\{u\leq X\}}|\mathbf{y}_i].$$

The proof is given in the Appendix. In this proof and in the following, we shall use the notation

$$L^{(i)}_{\hat{\theta}_n}(\theta) = E_{\hat{\theta}_n}[\ln l(\mathbf{Y}, Z; \theta)|\mathbf{y}_i],$$

and refer to $L_{n,\hat{\theta}_n}(\theta) = \sum_{i=1}^{n} L^{(i)}_{\hat{\theta}_n}(\theta)$ as the EM-loglikelihood.

## 4. Large sample properties.

4.1. *Consistency.* Since we are interested in almost sure (a.s.) consistency, we work with fixed realizations of the data which are assumed to lie in a set of probability one. Let $\|\cdot\|_\infty$ denote the supremum norm on $[0,\tau]$ and $\|\cdot\|$ denote the Euclidean norm.

THEOREM 1. *Under conditions* C1–C7, *the NPML estimator* $\hat{\theta}_n = (\hat{\alpha}_n, \hat{\beta}_n, \hat{\Lambda}_n)$ *is consistent:* $\|\hat{\alpha}_n - \alpha_0\|, |\hat{\beta}_n - \beta_0|$ *and* $\|\hat{\Lambda}_n - \Lambda_0\|_\infty$ *converge a.s. to zero as* $n \longrightarrow 0$.

PROOF. In the following, it will be convenient to denote $(\alpha, \beta)$ by $\gamma$. Our proof follows Murphy's [23] proof of a.s. consistency in the frailty model. The plan for proving consistency is as follows. We first show that the set $\{\hat{\theta}_n = (\hat{\gamma}_n, \hat{\Lambda}_n), n \in \mathbb{N}\}$ is relatively compact. Using the proposition on identifiability, we then show that its closure reduces to the single element $\theta_0 = (\gamma_0, \Lambda_0)$.

We first show that $(\hat{\Lambda}_n)_{n \in \mathbb{N}}$ stays bounded as $n \longrightarrow \infty$. We note from (A.2) in the Appendix that $\hat{\Lambda}_n(t) \geq 0$ for all $t \in [0,\tau]$ and that

$$\hat{\Lambda}_n(\tau) \leq \frac{\sum_{i=1}^{n}\Delta_i\mathbb{1}_{\{X_i\leq \tau\}}}{m\sum_{i=1}^{n}\mathbb{1}_{\{\tau\leq X_i\}}},$$

where $m = \min_{B,i,k} e^{\beta z_i(x_k)}$. Noting that there exists a constant $l$ such that $1/[\frac{1}{n}\sum_{i=1}^{n}\mathbb{1}_{\{\tau\leq X_i\}}] \leq 1/P(X \geq \tau) + l$ as $n \longrightarrow \infty$, it follows that $\hat{\Lambda}_n(\tau)$ does not diverge to infinity.

Let $\phi(n)$ be an arbitrary subsequence of $(n)$. From the Bolzanno–Weierstrass theorem, $(\hat{\gamma}_{\phi(n)})_{n\in\mathbb{N}}$ being a bounded sequence of $\mathbb{R}^{p+1}$ has a convergent subsequence $(\hat{\gamma}_{\varphi(\phi(n))})_{n\in\mathbb{N}}$ which converges to some $\gamma^*$. Since $\hat{\Lambda}_n$ is not allowed to diverge, the Helly–Bray lemma can be used to prove the existence of a subsequence $(\hat{\Lambda}_{\eta(\varphi(\phi(n)))})_{n\in\mathbb{N}}$ of $(\hat{\Lambda}_{\varphi(\phi(n))})_{n\in\mathbb{N}}$ which converges pointwise to some $\Lambda^*$. Since every subsequence of a convergent sequence in $\mathbb{R}^{p+1}$ must converge to the same limit, $(\hat{\gamma}_{\eta(\varphi(\phi(n)))})_{n\in\mathbb{N}}$ must converge to $\gamma^*$.

Hence, for any given subsequence $\hat{\theta}_{\phi(n)}$, we can find a further subsequence $\hat{\theta}_{\eta(\varphi(\phi(n)))}$ which converges to some $\theta^* = (\gamma^*, \Lambda^*)$. We now show that $\hat{\Lambda}_{\eta(\varphi(\phi(n)))}$ converges uniformly to $\Lambda^*$. In the following, we shall use the following notation for the sake of clarity of formulas: $g(n) = \eta(\varphi(\phi(n)))$.

We shall use in the sequel the Helly–Bray lemma and the result that the class of all functions $f : [0, \tau] \longrightarrow \mathbb{R}$ that are uniformly bounded and of variation bounded is Glivenko–Cantelli [33].

We first define the intermediate quantity $\bar{\Lambda}_n$ by $\bar{\Lambda}_n(t) = \int_0^t \frac{dH_n(u)}{W_n(u,\theta_0)}$, which will help mediate between $\hat{\Lambda}_n$ and $\Lambda_0$. By Glivenko–Cantelli, $(H_n(u))_{n\in\mathbb{N}}$ converges uniformly on $[0, \tau]$ to $H(u) = E_{\theta_0}[\Delta \mathbb{1}_{\{X \leq u\}}]$. Note that $\Lambda_0(t) = \int_0^t \frac{dH(u)}{W(u,\theta_0)}$, where $W(u, \theta_0) = E_{\theta_0}[e^{\beta_0 Z(u)} \mathbb{1}_{\{u \leq X\}}]$.

The functions $u \longmapsto E_{\theta_0}[e^{\beta_0 Z(u)} \mathbb{1}_{\{u \leq X\}} | \mathbf{y}]$ are uniformly bounded and of variation bounded. Hence, $(W_n(u, \theta_0))_{n\in\mathbb{N}}$ converges uniformly to $W(u, \theta_0) = E_{\theta_0}[e^{\beta_0 Z(u)} \mathbb{1}_{\{u \leq X\}}]$, which is bounded away from 0 by condition C6. Hence, $(1/W_n(u, \theta_0))_{n\in\mathbb{N}}$ converges to $1/W(u, \theta_0)$ uniformly on $[0, \tau]$.

Applying the Helly–Bray lemma gives that both $\|\bar{\Lambda}_n - \Lambda_0\|_\infty$ and $\|\hat{\Lambda}_{g(n)} - \Lambda^*\|_\infty$ converge almost surely to 0. This establishes the relative compactness of $\{\hat{\theta}_n, n \in \mathbb{N}\}$.

We now show that every subsequence $(\hat{\theta}_{g(n)})_{n\in\mathbb{N}}$ must converge to the true value $\theta_0 = (\gamma_0, \Lambda_0)$. Since $\hat{\theta}_{g(n)}$ maximizes the loglikelihood,

$$\frac{1}{g(n)} \sum_{i=1}^{g(n)} [\ln L^{(i)}(\hat{\theta}_{g(n)}) - \ln L^{(i)}(\gamma_0, \bar{\Lambda}_{g(n)})] \geq 0.$$

Note that, for all $g(n)$, as $m \longrightarrow \infty$,

$$\frac{1}{m} \sum_{i=1}^{m} [\ln L^{(i)}(\hat{\theta}_{g(n)}) - \ln L^{(i)}(\gamma_0, \bar{\Lambda}_{g(n)})]$$

$$\longrightarrow E_{\theta_0}[\ln L(\hat{\theta}_{g(n)}) - \ln L(\gamma_0, \bar{\Lambda}_{g(n)})] \qquad \text{a.s.}$$

It follows that

(4.4) $\qquad E_{\theta_0}[\ln L(\hat{\theta}_{g(n)}) - \ln L(\gamma_0, \bar{\Lambda}_{g(n)})] \geq -o(1).$

We have

$$\ln L(\hat{\theta}_{g(n)}) - \ln L(\gamma_0, \bar{\Lambda}_{g(n)}) \longrightarrow \ln L(\theta^*) - \ln L(\theta_0) \qquad \text{a.s.}$$

By Lebesgue's theorem,

$$E_{\theta_0}[\ln(L(\hat{\theta}_{g(n)})/L(\gamma_0, \bar{\Lambda}_{g(n)}))] \longrightarrow E_{\theta_0}[\ln(L(\theta^*)/L(\theta_0))] \qquad \text{a.s.}$$

From (4.4), $E_{\theta_0}[\ln(L(\theta^*)/L(\theta_0))] \geq 0$. This quantity cannot be treated directly as a Kullback–Leibler distance, since (2.2) is not a likelihood in the traditional sense. Due to the missing observation, it may even not be viewed as a generalized likelihood in the sense of Jacobsen [18]. However, it can be shown that $E_{\theta_0}[\ln(L(\theta^*)/L(\theta_0))] \leq 0$, and, moreover, that it is equal to zero if and only if $L(\theta^*) = L(\theta_0)$ a.e. This in turn implies $\theta^* = \theta_0$ by Proposition 1.
□

4.2. *Asymptotic normality.* In order to establish the asymptotic distribution of the proposed estimators $\hat{\theta}_n$, we follow the function analytic approach described by Murphy [24] to derive asymptotic theory for the frailty model. To calculate the score equations, instead of differentiating $L_{n,\theta}(\theta)$ with respect to $\alpha, \beta$ and the jump sizes of the cumulative baseline hazard function, we consider one-dimensional submodels through the estimators and we differentiate at the estimators. That is, we set $\theta_t = (\alpha_t, \beta_t, \Lambda_t)$, $\alpha_t = \alpha + th_1$, $\beta_t = \beta + th_2$ and $\Lambda_t(\cdot) = \int_0^\cdot (1 + th_3(u)) \, d\Lambda(u)$, where $h_1$ is a $p$-dimensional vector, $h_2 \in \mathbb{R}$, and $h_3$ is a function of bounded variation defined on $[0, \tau]$.

More precisely, let the class of $(h_1, h_2, h_3)$ be the space $H = \{h = (h_1, h_2, h_3) | h_1$ is a $p$-vector, $h_2 \in \mathbb{R}, h_3$ is a bounded function of bounded variation on $[0, \tau]\}$.

The following proposition gives the form of the empirical score $S_{n,\theta}$. Its proof is given in the Appendix. Define first

$$S_{n,\tilde{\theta},1}(\theta) = n^{-1} \sum_{i=1}^n E_{\tilde{\theta}}\left[\frac{\partial}{\partial \alpha} \ln f(Z_0, \ldots, Z_{a_X}, Z; \alpha) | \mathbf{y}_i\right],$$

$$S_{n,\tilde{\theta},2}(\theta) = n^{-1} \sum_{i=1}^n E_{\tilde{\theta}}\left[\Delta Z - \int_0^X Z(u) e^{\beta Z(u)} \, d\Lambda(u) | \mathbf{y}_i\right],$$

$$S_{n,\tilde{\theta},3}(\theta)(h_3) = n^{-1} \sum_{i=1}^n \left[\delta_i h_3(x_i) - E_{\tilde{\theta}}\left[\int_0^X h_3(u) e^{\beta Z(u)} \, d\Lambda(u) | \mathbf{y}_i\right]\right],$$

for some value $\tilde{\theta}$ of $\theta$, and the notation $S_{n,\tilde{\theta},12}^T(\theta) = (S_{n,\tilde{\theta},1}^T, S_{n,\tilde{\theta},2})(\theta)$, $S_{n,\tilde{\theta},12}(\theta) = n^{-1} \sum_{i=1}^n S_{\tilde{\theta},12}^{(i)}(\theta)$ and $h_{12}^T = (h_1^T, h_2)$.

PROPOSITION 4. *The empirical score can be written as*
$$S_{n,\tilde{\theta}}(\theta)(h) = h_{12}^T S_{n,\tilde{\theta},12}(\theta) + S_{n,\tilde{\theta},3}(\theta)(h_3).$$

In the sequel, letting $s_{\tilde{\theta}}(\mathbf{y}, \theta)(h) = \frac{\partial}{\partial t} L_{\tilde{\theta}}(\theta_t)|_{t=0}$, we shall write the empirical and expected scores as

$$S_{n,\tilde{\theta}}(\theta)(h) = \frac{1}{n} \sum_{i=1}^n s_{\tilde{\theta}}(\mathbf{y}_i, \theta)(h), \qquad S_{\tilde{\theta}}(\theta)(h) = E_{\theta_0}[s_{\tilde{\theta}}(\mathbf{Y}, \theta)(h)].$$

We define the following norm on $H$: if $h \in H$, let $\|h\|_H = \|h_1\| + |h_2| + \|h_3\|_v$, where $\|\cdot\|$ is the Euclidean norm and $\|h_3\|_v$ is the absolute value of $h_3(0)$ plus the total variation of $h_3$ on the interval $[0,\tau]$. We further define $H_p = \{h \in H, \|h\|_H \leq p\}$, $H_\infty = \{h \in H, \|h\|_H < \infty\}$ and $BV_p$ to be the space of real-valued functions on $[0,\tau]$ bounded by $p$ and of variation bounded by $p$.

Define $\theta(h) = (\alpha, \beta, \Lambda)(h) = h_1^T \alpha + h_2 \beta + \int_0^\tau h_3(u) \, d\Lambda(u)$. Then we can consider the parameter $\theta$ as a functional on $H_p$, and the parameter space $\Theta$ as a subset of $l^\infty(H_p)$, the space of bounded real-valued functions on $H_p$. We define on $l^\infty(H_p)$ the norm $\|U\|_p = \sup_{h \in H_p} |U(h)|$. For any finite $p$, the score function $S_{n,\theta}$ is a map from $\Theta$ to $l^\infty(H_p)$.

We obtain the following result:

THEOREM 2. *Let $0 < p < \infty$. Under assumptions C1–C7, the sequence $\sqrt{n}(\hat{\alpha}_n - \alpha_0, \hat{\beta}_n - \beta_0, \hat{\Lambda}_n - \Lambda_0)$ weakly converges in $l^\infty(H_p)$ to a centered Gaussian process $G$ with covariance process*

$$\mathrm{cov}[G(g), G(g^*)] = \int_0^\tau g_3(u) \sigma_{3,\theta_0}^{-1}(g^*)(u) \, d\Lambda_0(u) + \sigma_{2,\theta_0}^{-1}(g^*) g_2 + \sigma_{1,\theta_0}^{-1}(g^*)^T g_1,$$

*where $\sigma_{\theta_0}^{-1} = (\sigma_{1,\theta_0}^{-1}, \sigma_{2,\theta_0}^{-1}, \sigma_{3,\theta_0}^{-1})$ is the inverse of the continuously invertible linear operator $\sigma_{\theta_0} = (\sigma_{1,\theta_0}, \sigma_{2,\theta_0}, \sigma_{3,\theta_0})$ from $H_\infty$ to $H_\infty$, defined by*

$$\sigma_{1,\theta_0}(h) = -E_{\theta_0}\left[\frac{\partial^2}{\partial \alpha \, \partial \alpha^T} \ln f(Z_0, \ldots, Z_{a_X}, Z; \alpha_0)\right] h_1,$$

$$\sigma_{2,\theta_0}(h) = E_{\theta_0}\left[\int_0^X Z(u) e^{\beta_0 Z(u)} (Z(u) h_2 + h_3(u)) \, d\Lambda_0(u)\right],$$

$$\sigma_{3,\theta_0}(h)(u) = E_{\theta_0}\left[(Z(u) h_2 + h_3(u)) e^{\beta_0 Z(u)} \mathbb{1}_{\{u \leq X\}}\right].$$

PROOF. The proof is based on a theorem by van der Vaart and Wellner [33], which is stated as Lemma A.1 in the Appendix. In the following lemmas, we verify that the conditions stated in this theorem are satisfied by our estimator. Some additional technical lemmas are given in the Appendix, in order to keep attention on the main steps of the demonstration.

We first establish Fréchet differentiability of the map $\theta \longmapsto S_{\theta_0}(\theta)$ at $\theta_0$. Let us define the operator $\sigma_\theta$ from $H_\infty$ to $H_\infty$ by $\sigma_\theta(h) = (\sigma_{1,\theta}(h), \sigma_{2,\theta}(h), \sigma_{3,\theta}(h))$, where

$$\sigma_{1,\theta}(h) = -E_{\theta_0}\left[E_\theta\left[\frac{\partial^2}{\partial \alpha \, \partial \alpha^T} \ln f(Z_0, \ldots, Z_{a_X}, Z; \alpha) | \mathbf{y}\right]\right] h_1,$$

(4.5) $$\sigma_{2,\theta}(h) = E_{\theta_0}\left[E_\theta\left[\int_0^X Z(u) e^{\beta Z(u)} (Z(u) h_2 + h_3(u)) \, d\Lambda(u) | \mathbf{y}\right]\right],$$

$$\sigma_{3,\theta}(h)(u) = E_{\theta_0}[E_\theta[(Z(u) h_2 + h_3(u)) e^{\beta Z(u)} \mathbb{1}_{\{u \leq X\}} | \mathbf{y}]].$$

LEMMA 1. *For any finite $p$, the following holds: there exists a continuous linear operator $\dot{S}_{\theta_0}(\theta_0): \operatorname{lin}\Theta \longrightarrow l^\infty(H_p)$ such that $\|S_{\theta_0}(\theta) - S_{\theta_0}(\theta_0) - \dot{S}_{\theta_0}(\theta_0)(\theta - \theta_0)\|_p = o_P(\|\theta - \theta_0\|_p)$ as $\|\theta - \theta_0\|_p \longrightarrow 0$. The form of $\dot{S}_{\theta_0}(\theta_0)$ is as follows:*

$$\dot{S}_{\theta_0}(\theta_0)(\theta)(h) = -\int_0^\tau \sigma_{3,\theta_0}(h)(u)\,d\Lambda(u) - \beta\sigma_{2,\theta_0}(h) - \alpha^T \sigma_{1,\theta_0}(h).$$

PROOF. To establish this, we use the following characterization of Fréchet differentiability (see [2], page 454). Let $T$ be a function from a normed linear space $\mathbf{X}$ to another normed linear space $\mathbf{Y}$. Let $\mathcal{S}$ be the set of all bounded subsets of $\mathbf{X}$. $T$ is Fréchet differentiable at $x$ with derivative $\dot{T}_x$ if, for all $S \in \mathcal{S}$,

$$\frac{T(x+\varepsilon s) - T(x) - \dot{T}_x(\varepsilon s)}{\varepsilon} \longrightarrow 0 \quad \text{as } \varepsilon \longrightarrow 0 \text{ uniformly in } s \in S.$$

We first calculate the derivative $D_\theta S_{\theta_0}(\theta_0)$ given by

$$D_\theta S_{\theta_0}(\theta_0) = \frac{\partial}{\partial t} S_{\theta_0}(\theta_0 + t\theta)|_{t=0},$$

where $\theta_0 + t\theta = (\alpha_0 + t\alpha, \beta_0 + t\beta, \Lambda_0(\cdot) + t\Lambda(\cdot))$:

$$\frac{\partial}{\partial t} S_{\theta_0}(\theta_0 + t\theta)(h)$$
$$= \frac{\partial}{\partial t} E_{\theta_0}\bigg[\Delta h_3(X) - \int_0^X h_3(u) e^{[\beta_0 + t\beta]Z(u)}(d\Lambda_0(u) + t\,d\Lambda(u))$$
$$+ h_2\Delta Z - \int_0^X h_2 Z(u) e^{[\beta_0 + t\beta]Z(u)}(d\Lambda_0(u) + t\,d\Lambda(u))$$
$$+ h_1^T \frac{\partial}{\partial t} \ln f(Z_0, \ldots, Z_{a_X}, Z; \alpha_0 + t\alpha)\bigg].$$

The expression for $D_\theta S_{\theta_0}(\theta_0)(h)$ immediately follows and, using a first-order Taylor expansion of $\exp([\beta_0 + \varepsilon\beta]Z(u))$ around $\exp(\beta_0 Z(u))$, it is fairly straightforward to see that

$$S_{\theta_0}(\theta_0 + \varepsilon\theta)(h) - S_{\theta_0}(\theta_0)(h) - D_{\varepsilon\theta} S_{\theta_0}(\theta_0)(h) = o(\varepsilon).$$

Now, as $\varepsilon \longrightarrow 0$, $\frac{\|S_{\theta_0}(\theta_0 + \varepsilon\theta) - S_{\theta_0}(\theta_0) - D_{\varepsilon\theta} S_{\theta_0}(\theta_0)\|_p}{\varepsilon}$ converges to 0 uniformly in $\theta$ ranging over any element of the class of bounded subsets of $\operatorname{lin}\Theta$, where the notation "*lin*" before a set denotes the set of all finite linear combinations of elements of this set.

It follows that $S_{\theta_0}$ is Fréchet differentiable at $\theta_0$ and that the Fréchet derivative $\dot{S}_{\theta_0}(\theta_0)(\theta)$ is given by $\dot{S}_{\theta_0}(\theta_0)(\theta) = D_\theta S_{\theta_0}(\theta_0)$. □

We now consider the asymptotic distribution of the score function.

LEMMA 2. *For any finite $p$, the following holds: let $\mathcal{G}$ be a tight Gaussian process on $l^\infty(H_p)$ with covariance*

$$\operatorname{cov}(\mathcal{G}(h), \mathcal{G}(h^*)) = \int_0^\tau h_3(u)\sigma_{3,\theta_0}(h^*)(u)\,d\Lambda_0(u) + h_2\sigma_{2,\theta_0}(h^*) + h_1^T\sigma_{1,\theta_0}(h^*).$$

*Then*

$$\sqrt{n}(S_{n,\hat{\theta}_n}(\theta_0) - S_{\theta_0}(\theta_0)) \Longrightarrow \mathcal{G}.$$

PROOF. Note that $\sqrt{n}(S_{n,\hat{\theta}_n}(\theta_0) - S_{\theta_0}(\theta_0))(h)$ can be written as

$$\frac{1}{\sqrt{n}} \sum_{i=1}^n \left[ h_{12}^T S_{\hat{\theta}_n,12}^{(i)}(\theta_0) + \delta_i h_3(x_i) - \int_0^{x_i} h_3(u) E_{\hat{\theta}_n}[e^{\beta_0 Z(u)}|\mathbf{y}_i]\,d\Lambda_0(u) \right].$$

Note that $\{h_{12}^T S_{\hat{\theta}_n,12}(\theta_0) : h_1 \in \mathbb{R}^p, \|h_1\| \leq p, h_2 \in \mathbb{R}, |h_2| \leq p\}$ is bounded Donsker. The class $\{\delta h_3(x), h_3 \in BV_p\}$ is Donsker (this follows from the fact that the class of real-valued functions on $[0,\tau]$ that are uniformly bounded and are of variation bounded is Donsker). The class $\{\int_0^x h_3(u) E_{\hat{\theta}_n} \times [e^{\beta_0 Z(u)}|\mathbf{y}]\,d\Lambda_0(u) : h_3 \in BV_p\}$ is a bounded Donsker class. Then the class $\{h_{12}^T S_{\hat{\theta}_n,12}(\hat{\theta}_n) + \delta h_3(x) - \int_0^x h_3(u) E_{\hat{\theta}_n}[e^{\beta_0 Z(u)}|\mathbf{y}]\,d\Lambda_0(u) : h \in H_p\}$ is Donsker since the sum of bounded Donsker classes is Donsker. It follows that $\sqrt{n} S_{n,\hat{\theta}_n}(\hat{\theta}_n)$ converges in distribution to a zero mean tight Gaussian process $\mathcal{G}$ in $l^\infty(H_p)$.

The asymptotic distribution of the score $\sqrt{n}(S_{n,\hat{\theta}_n}(\theta_0) - S_{\theta_0}(\theta_0))$ is that of a tight Gaussian process $\mathcal{G}$ in $l^\infty(H_p)$ whose variance $\operatorname{var}(\mathcal{G}(h))$ is calculated as

$$-\frac{\partial}{\partial s} E_{\theta_0}[s_{\theta_0}(\mathbf{Y}, \theta_{0,s})(h)|_{s=0}] \qquad \text{which is } -\dot{S}_{\theta_0}(\theta_0)(h)(h).$$

The covariance of $\mathcal{G}$ is calculated as

$$\operatorname{cov}(\mathcal{G}(h), \mathcal{G}(h^*)) = -E_{\theta_0}\left[\frac{\partial}{\partial s}\frac{\partial}{\partial t} L_{\theta_0}(\theta_{0,s,t})|_{s,t=0}\right]$$

$$= -\dot{S}_{\theta_0}(\theta_0)(h^*)(h).$$

Let $\theta_s = (\alpha_s, \beta_s, \Lambda_s)$ with $\alpha_s = \alpha + sh_1^*$, $\beta_s = \beta + sh_2^*$ and $\Lambda_s(\cdot) = \int_0^\cdot (1 + sh_3^*(u))\,d\Lambda(u)$. Then we can calculate $\frac{\partial}{\partial s} s_{\tilde{\theta}}(\mathbf{y}, \theta_s)(h)$ as

$$\frac{\partial}{\partial s} s_{\tilde{\theta}}(\mathbf{y}, \theta_s)(h)$$

$$= \frac{\partial}{\partial s}\left[\frac{\partial}{\partial t} L_{\tilde{\theta}}(\theta_{s,t})|_{t=0}\right],$$

$$= -E_{\tilde{\theta}}\left[\int_0^x h_3(u) e^{[\beta+sh_2^*]z(u)}[h_2^* z(u)(1+sh_3^*(u)) + h_3^*(u)]\,d\Lambda(u)\Big|\mathbf{y}\right]$$

$$- h_2 E_{\tilde{\theta}}\left[\int_0^x z(u)e^{[\beta+sh_2^*]z(u)}[h_2^* z(u)(1+sh_3^*(u)) + h_3^*(u)]\,d\Lambda(u)|\mathbf{y}\right]$$

$$+ h_1^T E_{\tilde{\theta}}\left[\frac{\partial^2 \ln f(z_0,\ldots,z_{a_x},Z;\alpha+sh_1^*)}{\partial s\,\partial\alpha}\bigg|\mathbf{y}\right].$$

Calculation of $\frac{\partial}{\partial s}s_{\tilde{\theta}}(\mathbf{y},\theta_s)(h)|_{s=0}$ is straightforward, and using notation defined above, it follows that

$$\mathrm{cov}(\mathcal{G}(h),\mathcal{G}(h^*)) = \int_0^\tau h_3(u)\sigma_{3,\theta_0}(h^*)(u)\,d\Lambda_0(u) + h_2\sigma_{2,\theta_0}(h^*) + h_1^T\sigma_{1,\theta_0}(h^*).$$

□

The approximation condition (A.3) in Lemma A.1 follows by the Donsker property of the class of functions $\{s_\theta(y,\theta)(h) - s_{\theta_0}(y,\theta_0)(h) : \|\theta - \theta_0\|_p < \varepsilon, h \in H_p\}$ for some $\varepsilon > 0$. Details can be found in [11].

We now consider continuous invertibility of $\dot{S}_{\theta_0}(\theta_0)$.

LEMMA 3. *For any finite $p$, $\dot{S}_{\theta_0}(\theta_0)$ is continuously invertible on its range.*

PROOF. Continuous invertibility of $\dot{S}_{\theta_0}(\theta_0)$ on its range for some $p$ is equivalent (see [2], page 418) to the fact that there exists some $l > 0$ such that

$$(4.6) \qquad \inf_{\theta\in\mathrm{lin}\,\Theta} \frac{\|\dot{S}_{\theta_0}(\theta_0)(\theta)\|_p}{\|\theta\|_p} > l.$$

To prove (4.6), we follow two steps. We first show that $\sigma_{\theta_0}$ is a continuously invertible operator from $H_\infty$ to $H_\infty$. We can achieve this by proving that $\sigma_{\theta_0}$ is one-to-one and that it can be written as the sum of a continuously invertible operator $\Sigma$ plus a compact operator ([32], page 424).

From Lemma A.3 in the Appendix, we know that $\sigma_{\theta_0}$ is one-to-one, hence, we want to show that $\sigma_{\theta_0}$ can be written as the sum of a continuously invertible operator $\Sigma$ and a compact linear operator. We define $\Sigma$ as

$$\Sigma(h) = \left(-E_{\theta_0}\left[\frac{\partial^2}{\partial\alpha\,\partial\alpha^T}\ln f(Z_0,\ldots,Z_{a_X},Z;\alpha_0)\right]h_1,\right.$$
$$\left.E_{\theta_0}\left[\int_0^X \{Z(u)\}^2 e^{\beta_0 Z(u)}\,d\Lambda_0(u)\right]h_2, E_{\theta_0}[e^{\beta_0 Z(u)}\mathbb{1}_{\{u\leq X\}}]h_3(u)\right).$$

From conditions C1–C7, it follows that $\Sigma^{-1}$ is a bounded linear operator and, hence, that $\Sigma$ is continuously invertible.

We now have to show that $\sigma_{\theta_0}(h) - \Sigma(h)$ is compact. Let $(h_n)_{n\in\mathbb{N}} = (h_{1n},h_{2n},h_{3n})_{n\in\mathbb{N}}$ be a sequence in $H_p$. By the definition of a compact operator [27], we must prove that there exists a convergent subsequence of $\sigma_{\theta_0}(h_n) - \Sigma(h_n)$.

Since $h_{3n}$ is of bounded variation, we can write $h_{3n}$ as the difference of bounded increasing functions $h_{3n}^{(1)}$ and $h_{3n}^{(2)}$. From Helly's theorem, there exists a subsequence $(h_{3\phi(n)}^{(1)})$ of $(h_{3n}^{(1)})$ which converges pointwise to some $h_3^{(1)*}$. There also exists a subsequence $(h_{3\eta(\phi(n))}^{(2)})$ of $(h_{3\phi(n)}^{(2)})$ which converges pointwise to some $h_3^{(2)*}$. Finally, $(h_{3\eta(\phi(n))}^{(1)}, h_{3\eta(\phi(n))}^{(2)})$ converges pointwise to $h_3^* = (h_3^{(1)*}, h_3^{(2)*})$. Using the same argument and the Bolzanno–Weierstrass theorem, we can find a subsequence of $(h_n)_{n\in\mathbb{N}}$ [let us denote it by $(h_{g(n)})_{n\in\mathbb{N}}$ for notational simplicity] that converges to $h^* = (h_1^*, h_2^*, h_3^*)$.

We must prove that $\sigma_{\theta_0}(h_{g(n)}) - \Sigma(h_{g(n)})$ converges to $\sigma_{\theta_0}(h^*) - \Sigma(h^*)$ in $H_p$ for all $p$. Note that $\sigma_{\theta_0}(h) - \Sigma(h)$ is equal to

$$\left(0, E_{\theta_0}\left[\int_0^X Z(u)e^{\beta_0 Z(u)}h_3(u)\,d\Lambda_0(u)\right], E_{\theta_0}[h_2 Z(u)e^{\beta_0 Z(u)}\mathbb{1}_{\{u\leq X\}}]\right).$$

Now $\|\sigma_{\theta_0}(h_{g(n)}) - \Sigma(h_{g(n)}) - \sigma_{\theta_0}(h^*) + \Sigma(h^*)\|_H$ is equal to

$$\left|E_{\theta_0}\left[\int_0^X Z(u)e^{\beta_0 Z(u)}(h_{3g(n)} - h_3^*)(u)\,d\Lambda_0(u)\right]\right|$$
$$+ \|E_{\theta_0}[(h_{2g(n)} - h_2^*)Z(u)e^{\beta_0 Z(u)}\mathbb{1}_{\{u\leq X\}}]\|_v,$$

which, under conditions C1–C7, is bounded above by

$$ce^{bc} \cdot \int_0^\tau |(h_{3g(n)} - h_3^*)(u)|\,d\Lambda_0(u) + ce^{bc}(2+c) \cdot |h_{2g(n)} - h_2^*|,$$

where $b$ is such that $|\beta| < b$. From the dominated convergence theorem, the first term converges to zero and the overall bound converges to zero. It follows that $\sigma_{\theta_0}(h) - \Sigma(h)$ is a compact operator for all $p$.

We have then proved that $\sigma_{\theta_0}$ is a continuously invertible operator. This means that, for all $p > 0$, there exists a $q > 0$ such that $\sigma_{\theta_0}^{-1}(H_q) \subset H_p$. Hence, the LHS of (4.6) is bounded below by

$$\inf_{\theta\in\mathrm{lin}\Theta} \frac{\sup_{h\in\sigma_{\theta_0}^{-1}(H_q)} |\dot{S}_{\theta_0}(\theta_0)(\theta)(h)|}{\|\theta\|_p}$$
$$= \inf_{\theta\in\mathrm{lin}\,\Theta}\left[\sup_{h\in H_q}\left|\int_0^\tau \sigma_{3,\theta_0}(\sigma_{\theta_0}^{-1}(h))(u)\,d\Lambda(u)\right.\right.$$
$$\left.\left. + \beta\sigma_{2,\theta_0}(\sigma_{\theta_0}^{-1}(h)) + \alpha^T\sigma_{1,\theta_0}(\sigma_{\theta_0}^{-1}(h))\right|\right](\|\theta\|_p)^{-1}.$$

Recall that $\sigma_{\theta_0}$ is invertible, hence, $\sigma_{\theta_0}(\sigma_{\theta_0}^{-1}(h)) = h$ for all $h = (h_1, h_2, h_3) \in H_q$. Since $\sigma_{\theta_0}(\sigma_{\theta_0}^{-1}(h)) = (\sigma_{1,\theta_0}(\sigma_{\theta_0}^{-1}(h)), \sigma_{2,\theta_0}(\sigma_{\theta_0}^{-1}(h)), \sigma_{3,\theta_0}(\sigma_{\theta_0}^{-1}(h)))$, it follows that $\sigma_{i,\theta_0}(\sigma_{\theta_0}^{-1}(h)) = h_i$, $i = 1, 2, 3$. Hence, the above bound can be rewritten as

$$(4.7) \quad \inf_{\theta \in \lin \Theta} \frac{\sup_{h \in H_q} |\int_0^\tau h_3(u)\, d\Lambda(u) + h_2\beta + h_1^T \alpha|}{\|\theta\|_p} = \inf_{\theta \in \lin \Theta} \frac{\|\theta\|_q}{\|\theta\|_p}.$$

Now from Lemma A.2 in the Appendix, $\|\theta\|_q$ is greater than or equal to $q(\|\alpha\| \vee |\beta| \vee V_{[0,\tau]}(\Lambda))$ and $\|\theta\|_p$ is less than or equal to $3p(\|\alpha\| \vee |\beta| \vee V_{[0,\tau]}(\Lambda))$, where $V_{[0,\tau]}(f)$ denotes the total variation of a function $f$ on $[0, \tau]$. Hence, the RHS of (4.7) is greater than or equal to $\frac{q}{3p}$. It follows that $\dot{S}_{\theta_0}(\theta_0)$ is continuously invertible. □

Putting all these results together, it follows that, for all $h \in H_p$,

$$-\dot{S}_{\theta_0}(\theta_0)\sqrt{n}(\hat{\theta}_n - \theta_0)(h) = \sqrt{n}(S_{n,\hat{\theta}_n}(\theta_0) - S_{\theta_0}(\theta_0))(h) + o_P(1),$$

where

$$-\dot{S}_{\theta_0}(\theta_0)\sqrt{n}(\hat{\theta}_n - \theta_0)(h)$$
$$= \int_0^\tau \sigma_{3,\theta_0}(h)(u)\sqrt{n}\, d(\hat{\Lambda}_n - \Lambda_0)(u)$$
$$+ \sqrt{n}(\hat{\beta}_n - \beta_0)\sigma_{2,\theta_0}(h) + \sqrt{n}(\hat{\alpha}_n - \alpha_0)^T\sigma_{1,\theta_0}(h).$$

Consequently, $\sqrt{n}(\hat{\theta}_n - \theta_0)(h) \Longrightarrow -\dot{S}_{\theta_0}(\theta_0)^{-1}\mathcal{G}(h)$. We now want to identify $\dot{S}_{\theta_0}(\theta_0)^{-1}\mathcal{G}$. $\sigma_{\theta_0}$ is continuously invertible, hence, for all $q > 0$ there exists a $p > 0$ such that $\sigma_{\theta_0}^{-1}(g) \in H_q$ if $g \in H_p$. Let $h = \sigma_{\theta_0}^{-1}(g)$. Then

$$-\dot{S}_{\theta_0}(\theta_0)\sqrt{n}(\hat{\theta}_n - \theta_0)(h) = \int_0^\tau g_3(u)\sqrt{n}\, d(\hat{\Lambda}_n - \Lambda_0)(u)$$
$$(4.8) \qquad\qquad\qquad\qquad + \sqrt{n}(\hat{\beta}_n - \beta_0)g_2 + \sqrt{n}(\hat{\alpha}_n - \alpha_0)^T g_1$$

and

$$(4.9) \quad -\dot{S}_{\theta_0}(\theta_0)\sqrt{n}(\hat{\theta}_n - \theta_0)(h) = \sqrt{n}(S_{n,\theta_0}(\theta_0) - S_{\theta_0}(\theta_0))(\sigma_{\theta_0}^{-1}(g)) + o_P(1).$$

Note that the RHS of (4.8) is $(\sqrt{n}(\hat{\alpha}_n - \alpha_0), \sqrt{n}(\hat{\beta}_n - \beta_0), \sqrt{n}(\hat{\Lambda}_n - \Lambda_0))(g)$, which converges to $-\dot{S}_{\theta_0}(\theta_0)^{-1}\mathcal{G}(g)$. Note also that the RHS of (4.9) converges to $\mathcal{G}(\sigma_{\theta_0}^{-1}(g))$ [which has mean zero and variance $\int_0^\tau g_3(u)\sigma_{3,\theta_0}^{-1}(g)(u)\, d\Lambda_0(u) + \sigma_{2,\theta_0}^{-1}(g)g_2 + \sigma_{1,\theta_0}^{-1}(g)^T g_1$].

It then follows that $-\dot{S}_{\theta_0}(\theta_0)^{-1}\mathcal{G} = \mathcal{G}(\sigma_{\theta_0}^{-1})$. Hence, $(\sqrt{n}(\hat{\alpha}_n - \alpha_0), \sqrt{n}(\hat{\beta}_n - \beta_0), \sqrt{n}(\hat{\Lambda}_n - \Lambda_0))$ converges in $l^\infty(H_p)$ to a tight Gaussian process $G$ in $l^\infty(H_p)$ with mean zero and covariance process

$$\cov[G(g), G(g^*)] = \int_0^\tau g_3(u)\sigma_{3,\theta_0}^{-1}(g^*)(u)\, d\Lambda_0(u) + \sigma_{2,\theta_0}^{-1}(g^*)g_2 + \sigma_{1,\theta_0}^{-1}(g^*)^T g_1,$$

where $(\sigma_{1,\theta_0}^{-1}, \sigma_{2,\theta_0}^{-1}, \sigma_{3,\theta_0}^{-1})$ is the inverse of $\sigma_{\theta_0} = (\sigma_{1,\theta_0}, \sigma_{2,\theta_0}, \sigma_{3,\theta_0})$. $\square$

We now consider the problem of estimating the asymptotic variance of the NPML estimator. From Theorem 2, the asymptotic variance of

$$(\sqrt{n}(\hat{\alpha}_n - \alpha_0), \sqrt{n}(\hat{\beta}_n - \beta_0), \sqrt{n}(\hat{\Lambda}_n - \Lambda_0))(h)$$

is $\int_0^\tau h_3(u)\sigma_{3,\theta_0}^{-1}(h)(u)\,d\Lambda_0(u) + \sigma_{2,\theta_0}^{-1}(h)h_2 + \sigma_{1,\theta_0}^{-1}(h)^T h_1$. Using formulas (4.5), we propose to first estimate $\sigma_{\theta_0}$ by $\hat{\sigma}_{\hat{\theta}_n} = (\hat{\sigma}_{1,\hat{\theta}_n}, \hat{\sigma}_{2,\hat{\theta}_n}, \hat{\sigma}_{3,\hat{\theta}_n})$, where

$$\hat{\sigma}_{1,\hat{\theta}_n}(h) = -\frac{1}{n}\sum_{i=1}^n E_{\hat{\theta}_n}\left[\frac{\partial^2}{\partial\alpha\,\partial\alpha^T}\ln f(Z_0,\ldots,Z_{a_X},Z;\hat{\alpha}_n)|\mathbf{y}_i\right]h_1,$$

$$\hat{\sigma}_{2,\hat{\theta}_n}(h) = \frac{1}{n}\sum_{i=1}^n E_{\hat{\theta}_n}\left[\int_0^X Z(u)e^{\hat{\beta}_n Z(u)}[h_2 Z(u) + h_3(u)]\,d\hat{\Lambda}_n(u)|\mathbf{y}_i\right],$$

$$\hat{\sigma}_{3,\hat{\theta}_n}(h)(u) = \frac{1}{n}\sum_{i=1}^n E_{\hat{\theta}_n}[[h_2 Z(u) + h_3(u)]e^{\hat{\beta}_n Z(u)}\mathbb{1}_{\{u \leq X\}}|\mathbf{y}_i].$$

We then propose to estimate the asymptotic variance by

$$(4.10) \quad \int_0^\tau h_3(u)\hat{\sigma}_{3,\hat{\theta}_n}^{-1}(h)(u)\,d\hat{\Lambda}_n(u) + \hat{\sigma}_{2,\hat{\theta}_n}^{-1}(h)h_2 + \hat{\sigma}_{1,\hat{\theta}_n}^{-1}(h)^T h_1.$$

Using the same arguments as in the proof of Lemma 2, we are able to show that the functions under $\sum$ in $\hat{\sigma}_{1,\hat{\theta}_n}, \hat{\sigma}_{2,\hat{\theta}_n}, \hat{\sigma}_{3,\hat{\theta}_n}$ form Donsker classes (for $h \in H_p$, $p > 0$), hence $\sup_{h \in H_p}\|\hat{\sigma}_{\hat{\theta}_n}(h) - \sigma_{\theta_0}(h)\|_H \longrightarrow 0$.

$\hat{\sigma}_{\hat{\theta}_n}$ is continuously invertible, hence, for all $H_p \subset H_\infty$ there exists $H_q \subset H_\infty$ such that $\hat{\sigma}_{\hat{\theta}_n}^{-1}(H_q) \subset H_p$, and for all $g \in H_q$ there exists $h \in H_p$ such that $h = \hat{\sigma}_{\hat{\theta}_n}^{-1}(g)$.

Then

$$\|\hat{\sigma}_{\hat{\theta}_n}^{-1}(g) - \sigma_{\theta_0}^{-1}(g)\|_H = \|\sigma_{\theta_0}^{-1}(\sigma_{\theta_0}(h)) - \sigma_{\theta_0}^{-1}(\hat{\sigma}_{\hat{\theta}_n}(h))\|_H$$

$$\leq \sup_{h \in H_q}\frac{\|\sigma_{\theta_0}^{-1}(h)\|_H}{\|h\|_H}\sup_{h \in H_p}\|\sigma_{\theta_0}(h) - \hat{\sigma}_{\hat{\theta}_n}(h)\|_H.$$

It follows that $\sup_{g \in H_q}\|\hat{\sigma}_{\hat{\theta}_n}^{-1}(g) - \sigma_{\theta_0}^{-1}(g)\|_H \longrightarrow 0$ and that the sequence of estimators (4.10) converges to the limit $\int_0^\tau h_3(u)\sigma_{3,\theta_0}^{-1}(h)(u)\,d\Lambda_0(u) + \sigma_{2,\theta_0}^{-1}(h)h_2 + \sigma_{1,\theta_0}^{-1}(h)^T h_1$.

In the above framework, specific choices of $h$ allow one to estimate the asymptotic variance of any particular estimator. For example, by setting $h_\beta = (0, 1, 0)$, one may obtain the following convergent estimator of the asymptotic variance of $\sqrt{n}(\hat{\beta}_n - \beta_0)$:

$$\hat{\sigma}^{-1}_{2,\hat{\theta}_n}(h_\beta) = \left[\frac{1}{n}\sum_{i=1}^{n}\sum_{x_k \leq x_i} E_{\hat{\theta}_n}[\{Z(x_k)\}^2 e^{\hat{\beta}_n Z(x_k)}|\mathbf{y}_i]\Delta\hat{\Lambda}_{n,k}\right]^{-1}.$$

**5. Discussion.** In this paper we have described a joint modeling approach for estimation in the Cox model with missing values of a time-dependent covariate. We have used nonparametric maximum likelihood estimation. Using the theory of empirical processes and techniques developed by Murphy [23, 24] and van der Vaart and Wellner [33], we have shown that the proposed estimators are consistent and asymptotically normal. Moreover, we have proposed a consistent estimator of the asymptotic variance.

An alternative widely used approach to the modeling of longitudinal data in joint models assumes a normal random effects model for the repeated measurements (see among others Henderson, Diggle and Dobson [17], Tsiatis, DeGruttola and Wulfsohn [31] and Wulfsohn and Tsiatis [34]). Tsiatis and Davidian [29] estimate parameters in a joint model without requiring any distributional assumption on the random effects, and an informal proof of large-sample properties of estimators in the proportional hazards model is given. However, a formal theoretical justification of asymptotic properties for the maximum likelihood estimator in joint models with random effects for longitudinal data is not yet available, and should be a subject for further work.

## APPENDIX

PROOF OF PROPOSITION 1. In the following, a.e. will stand for almost everywhere. Using (2.2), $\ln L(\theta) = \ln L(\theta')$ a.e. can be reexpressed as

$$\delta \ln\left(\frac{\lambda}{\lambda'}(x)\right) + \ln \int \exp\left[\delta\beta z - \int_t^x e^{\beta z(u)}\,d\Lambda(u)\right] f(z_0,\ldots,z_{a_x},z;\alpha)\,dz$$

$$(\text{A.1}) \quad -\ln \int \exp\left[\delta\beta' z - \int_t^x e^{\beta' z(u)}\,d\Lambda'(u)\right] f(z_0,\ldots,z_{a_x},z;\alpha')\,dz$$

$$= \int_0^t [e^{\beta z(u)}\,d\Lambda(u) - e^{\beta' z(u)}\,d\Lambda'(u)] \qquad \text{a.e.}$$

for $t < x$ and $x \in (0,\tau]$. The LHS of (A.1) depends on the path of the longitudinal covariate only through $\{z(u), t \leq u \leq x\}$ and $\{z_0,\ldots,z_{a_t}\}$. Hence, for given $\{z(u): t \leq u \leq x\}$ and $\{z_0,\ldots,z_{a_t}\}$, the RHS of (A.1) should yield the same value for two different paths $z(u)$ and $z^*(u)$ ($0 \leq u \leq t$) taking the same values $\{z_0,\ldots,z_{a_t}\}$ at $t_0,\ldots,t_{a_t}$. This can be expressed as

$$\int_0^t [e^{\beta z(u)}\,d\Lambda(u) - e^{\beta' z(u)}\,d\Lambda'(u)] = \int_0^t [e^{\beta z^*(u)}\,d\Lambda(u) - e^{\beta' z^*(u)}\,d\Lambda'(u)].$$

Letting $z(u) = \xi$ and $z^*(u) = \xi + h$ ($h > 0$) in $[0, t]$ except at $t_0, \ldots, t_{a_t}$ [where $z(u)$ and $z^*(u)$ take values $z_0, \ldots, z_{a_t}$], and eventually also at at most a countable number of time points in $[0, t]$, the following holds: $e^{(\beta - \beta')\xi} = \Lambda'(t)(1 - e^{\beta' h})/\Lambda(t)(1 - e^{\beta h})$. For a fixed $h$, the RHS of this expression is independent of $\xi$. It follows that $\beta = \beta'$ and then $\Lambda = \Lambda'$. Rewriting $\ln L(\theta) = \ln L(\theta')$ with $\beta = \beta'$ and $\Lambda = \Lambda'$ leads to $\alpha = \alpha'$. □

PROOF OF PROPOSITION 2. Recall that $\alpha$ and $\beta$ are interior points of some known compact sets $A$ and $B$. Suppose first that $\Delta \Lambda_{n,i} \leq U$ [$i = 1, \ldots, p(n)$] for some finite $U$. Since $L_n$ is a continuous function of $\alpha, \beta$ and of the jump sizes $\Delta \Lambda_{n,i}$ [$i = 1, \ldots, p(n)$] on the compact set $A \times B \times [0, U]^{p(n)}$, $L_n$ achieves its maximum on this set.

To show that a maximum of $L_n$ exists on $A \times B \times [0, \infty)^{p(n)}$, we show that there exists a finite $U$ such that, for all $\theta_U = (\alpha_U, \beta_U, \Delta \Lambda_{n,1,U}, \ldots, \Delta \Lambda_{n,p(n),U}) \in \{A \times B \times [0, \infty)^{p(n)}\} \setminus \{A \times B \times [0, U]^{p(n)}\}$, there exists a $\theta = (\alpha, \beta, \Delta \Lambda_{n,1}, \ldots, \Delta \Lambda_{n,p(n)}) \in A \times B \times [0, U]^{p(n)}$ such that $L_n(\theta) > L_n(\theta_U)$.

A proof by contradiction is adopted for this purpose. Assume that, for all $U$, there exists $\theta_U = (\alpha_U, \beta_U, \Delta \Lambda_{n,1,U}, \ldots, \Delta \Lambda_{n,p(n),U}) \in \{A \times B \times [0, \infty)^{p(n)}\} \setminus \{A \times B \times [0, U]^{p(n)}\}$ such that, for all $\theta = (\alpha, \beta, \Delta \Lambda_{n,1}, \ldots, \Delta \Lambda_{n,p(n)}) \in A \times B \times [0, U]^{p(n)}$, $L_n(\theta) < L_n(\theta_U)$. Under conditions C1–C7, it is easily seen that the likelihood $L_n(\theta)$ is bounded above by

$$\prod_{i=1}^{n} \left[ (M \Delta \Lambda_n(x_i))^{\delta_i} \exp\left( -m \sum_{k=1}^{p(n)} \Delta \Lambda_{n,k} \mathbb{1}_{\{x_k \leq x_i\}} \right) f(z_{i,0}, \ldots, z_{i,a_{x_i}}; \alpha) \right],$$

where $m = \min_{B,i,k} e^{\beta z_i(x_k)}$ and $M = \max_{B,i,k} e^{\beta z_i(x_k)}$.

If $\theta_U = (\alpha_U, \beta_U, \Delta \Lambda_{n,1,U}, \ldots, \Delta \Lambda_{n,p(n),U}) \in \{A \times B \times [0, \infty)^{p(n)}\} \setminus \{A \times B \times [0, U]^{p(n)}\}$, then there exists $j$ [$j \in \{1, \ldots, p(n)\}$] such that $\Delta \Lambda_{n,j,U} > U$. Hence, there exists at least one $i_U$ ($i_U \in \{1, \ldots, n\}$) such that $\sum_{k=1}^{p(n)} \Delta \Lambda_{n,k,U} \mathbb{1}_{\{x_k \leq x_{i_U}\}} > U$. Hence, $\sum_{k=1}^{p(n)} \Delta \Lambda_{n,k,U} \mathbb{1}_{\{x_k \leq x_{i_U}\}} \longrightarrow +\infty$ as $U \longrightarrow +\infty$. It follows that the upper bound of $L_n(\theta_U)$ [and, hence, $L_n(\theta_U)$] can be made as close to 0 as desired by increasing $U$. This is the desired contradiction. □

PROOF OF PROPOSITION 3. The maximizer $\hat{\theta}_n$ of $\sum_{j=1}^{n} \ln f_{\mathbf{Y}}(\mathbf{y}_j; \theta)$ over $\theta \in \Theta_n$ satisfies

$$\sum_{j=1}^{n} \frac{\partial}{\partial \theta} [E_{\hat{\theta}_n}[\ln f_{\mathbf{Y}, Z}(\mathbf{Y}, Z; \theta) | \mathbf{y}_j]]_{|\theta = \hat{\theta}_n} = 0.$$

This result can be obtained by using the same argument that Dempster, Laird and Rubin [9] used to derive the principle of the EM algorithm. Its proof is therefore omitted. Discarding from $E_{\hat{\theta}_n}[\ln f_{\mathbf{Y},Z}(\mathbf{Y}, Z; \theta)|\mathbf{y}_j]$ terms adhering to censoring does not influence maximization, hence we can write that

$$\sum_{j=1}^{n} \frac{\partial}{\partial \theta} L_{\hat{\theta}_n}^{(j)}(\theta)_{|\theta=\hat{\theta}_n} = 0.$$

Letting $\theta \in \Theta_n$, we note that $L_{\hat{\theta}_n}^{(j)}(\theta)$ is equal to

$$E_{\hat{\theta}_n}\left[\Delta\beta Z - \sum_{k=1}^{p(n)} \Delta\Lambda_{n,k} e^{\beta Z(X_k)} \mathbb{1}_{\{X_k \leq X\}} \right.$$
$$\left. + [\ln \Delta\Lambda_{n,j}]^\Delta + \ln f(Z_0, \ldots, Z_{a_X}, Z; \alpha)|\mathbf{y}_j\right].$$

Summing this expression over $j$ ($j = 1, \ldots, n$), deriving with respect to $\Delta\Lambda_{n,i}$ and solving the derivative to 0 gives $\Delta\hat{\Lambda}_{n,i} = 1/[\sum_{j=1}^{n} E_{\hat{\theta}_n}[e^{\hat{\beta}_n Z(X_i)} \mathbb{1}_{\{X_i \leq X\}}|\mathbf{y}_j]]$. The cumulative baseline hazard can be estimated by $\hat{\Lambda}_n(t) = \sum_{i=1}^{p(n)} \Delta\hat{\Lambda}_{n,i} \mathbb{1}_{\{X_i \leq t\}}$. Using $H_n$ and $W_n$ given in Proposition 3, this can further be written as

$$\text{(A.2)} \qquad \hat{\Lambda}_n(t) = \int_0^t \frac{dH_n(u)}{W_n(u; \hat{\theta}_n)}. \qquad \square$$

PROOF OF PROPOSITION 4. Letting $\tilde{\theta}$ be some value of $\theta$ and $\theta_t = (\alpha_t, \beta_t, \Lambda_t)$, $L_{\tilde{\theta}}(\theta_t)$ is equal to

$$L_{\tilde{\theta}}(\theta_t) = \delta \ln[(1 + th_3(x)) d\Lambda(x)] - E_{\tilde{\theta}}\left[\int_0^X e^{[\beta+th_2]Z(u)}(1 + th_3(u)) d\Lambda(u)|\mathbf{y}\right]$$
$$+ E_{\tilde{\theta}}[\Delta[\beta + th_2]Z|\mathbf{y}] + E_{\tilde{\theta}}[\ln f(Z_0, \ldots, Z_{a_X}, Z; \alpha + th_1)|\mathbf{y}].$$

Then

$$\frac{\partial}{\partial t} L_{\tilde{\theta}}(\theta_t) = \delta \frac{h_3(x)}{1 + th_3(x)} - E_{\tilde{\theta}}\left[\int_0^X h_3(u) e^{[\beta+th_2]Z(u)} d\Lambda(u)|\mathbf{y}\right]$$
$$+ h_2 E_{\tilde{\theta}}\left[\Delta Z - \int_0^X Z(u) e^{[\beta+th_2]Z(u)}(1 + th_3(u)) d\Lambda(u)|\mathbf{y}\right]$$
$$+ E_{\tilde{\theta}}\left[\frac{\partial}{\partial t} \ln f(Z_0, \ldots, Z_{a_X}, Z; \alpha + th_1)|\mathbf{y}\right].$$

Letting $t = 0$ in this derivative and using notation previously defined for $S_{n,\tilde{\theta},1}(\theta), S_{n,\tilde{\theta},2}(\theta)$ and $S_{n,\tilde{\theta},3}(\theta)$, Proposition 4 immediatly follows.

We verify that $S_{n,\hat{\theta}_n}(\hat{\theta}_n) = 0$. To see this, recall (proof of Proposition 3) that the maximizer $\hat{\theta}_n$ of $\sum_{i=1}^{n} \ln f_{\mathbf{Y}}(\mathbf{y}_i; \theta)$ over $\Theta_n$ satisfies

$$\frac{1}{n} \sum_{i=1}^{n} \frac{\partial}{\partial \theta} L_{\hat{\theta}_n}^{(i)}(\theta)|_{\theta=\hat{\theta}_n} = 0, \qquad \text{or, equivalently,} \quad \frac{1}{n} \sum_{i=1}^{n} s_{\hat{\theta}_n}(\mathbf{y}_i, \hat{\theta}_n)(h) = 0.$$

Hence, $S_{n,\hat{\theta}_n}(\hat{\theta}_n)(h) = 0$.

Note also that $S_{\theta_0}(\theta_0) = 0$. To see this, recall that, in Proposition 1 it was shown that the model is identifiable, that is, the function $t \longmapsto E_{\theta_0}[\ln f_{\mathbf{Y}}(\mathbf{Y};\theta_{0,t})]$ has a unique maximum at $t = 0$ [where $\theta_{0,t} = (\alpha_{0,t}, \beta_{0,t}, \Lambda_{0,t})$].

Since $E_{\theta_0}[\frac{\partial}{\partial t} \ln f_{\mathbf{Y}}(\mathbf{Y};\theta_{0,t})|_{t=0}] = 0$ and

$$\frac{\partial}{\partial t} \ln f_{\mathbf{Y}}(\mathbf{y};\theta_{0,t})|_{t=0} = \frac{\partial}{\partial t}[E_{\theta_0}[\ln f_{\mathbf{Y},Z}(\mathbf{Y}, Z;\theta_{0,t})|\mathbf{y}]]|_{t=0},$$

it follows that $E_{\theta_0}[\frac{\partial}{\partial t}[E_{\theta_0}[\ln f_{\mathbf{Y},Z}(\mathbf{Y}, Z;\theta_{0,t})|\mathbf{y}]]|_{t=0}] = 0$. By discarding terms adhering to censoring from $E_{\theta_0}[\ln f_{\mathbf{Y},Z}(\mathbf{Y}, Z;\theta_{0,t})|\mathbf{y}]$, it follows that $E_{\theta_0}[\frac{\partial}{\partial t} L_{\theta_0}(\theta_{0,t})|_{t=0}] = 0$. Hence $S_{\theta_0}(\theta_0)(h) = 0$. □

In the following lemma, we recall Theorem 3.3.1 of [33], which will serve as a basis for our proof of asymptotic normality.

LEMMA A.1 (Theorem 3.3.1 of [33]). *Let $S_n$ and $S$ be random maps and a fixed map, respectively, from $\Psi$ into a Banach space such that*

(A.3) $\quad \sqrt{n}(S_n - S)(\hat{\psi}_n) - \sqrt{n}(S_n - S)(\psi_0) = o_P(1 + \sqrt{n}\|\hat{\psi}_n - \psi_0\|),$

*and such that the sequence $\sqrt{n}(S_n - S)(\psi_0)$ converges in distribution to a tight random element $Z$. Let $\psi \longmapsto S(\psi)$ be Fréchet-differentiable at $\psi_0$ with a continuously invertible derivative $\dot{S}(\psi_0)$. If $S(\psi_0) = 0$ and $\hat{\psi}_n$ satisfies $S_n(\hat{\psi}_n) = o_P(n^{-1/2})$ and $\hat{\psi}_n - \psi_0 = o_P(1)$, then*

$$\sqrt{n}(\hat{\psi}_n - \psi_0) \Rightarrow -\dot{S}(\psi_0)^{-1} Z.$$

The following two technical lemmas will be useful for our proof.

LEMMA A.2. *For any finite $p$, the following holds: if $\theta \in l^{\infty}(H_p)$, then*

$$p(\|\alpha\| \vee |\beta| \vee V_{[0,\tau]}(\Lambda)) \leq \|\theta\|_p \leq 3p(\|\alpha\| \vee |\beta| \vee V_{[0,\tau]}(\Lambda)).$$

This lemma can easily be proved and its proof is therefore omitted. See [11] for details.

LEMMA A.3. *The operator $\sigma_{\theta_0}$ is such that $\mathrm{Ker}(\sigma_{\theta_0}) = \{0\}$.*

PROOF. Suppose that $\sigma_{\theta_0}(h) = 0$ for some $h = (h_1, h_2, h_3)$. Then $\sigma_{1,\theta_0}(h) = 0$. It follows that $h_1^T \sigma_{1,\theta_0}(h) = 0$. Since $-E_{\theta_0}[\frac{\partial^2}{\partial \alpha \partial \alpha^T} \ln f(Z_0, \ldots, Z; \alpha_0)]$ is positive definite (condition C6), it follows that $h_1 = 0$.

If we assume $\sigma_{\theta_0}(h) = 0$, then

$$\int_0^{\tau} h_3(u) \sigma_{3,\theta_0}(h)(u) \, d\Lambda_0(u) + h_2 \sigma_{2,\theta_0}(h) + h_1^T \sigma_{1,\theta_0}(h) = 0,$$

which we can rewrite as $E_{\theta_0}[(s_{\theta_0}(\mathbf{Y},\theta_0)(h))^2] = 0$. From this, it follows that $s_{\theta_0}(\mathbf{y},\theta_0)(h) = 0$ a.e. Using the fact that $h_1 = 0$, and acting similarly as in proof of identifiability, we can show that $h_2 = 0$.

At this stage, we get that $h_1 = 0$ and $h_2 = 0$. Let $h = (0, 0, h_3)$. Then

$$\sigma_{3,\theta_0}(h)(u) = h_3(u) E_{\theta_0}[e^{\beta_0 Z(u)} \mathbb{1}_{\{u \leq X\}}] = 0 \qquad \text{for all } u.$$

It follows from condition C6 that $h_3(u) = 0$ for all $u$. We conclude that $\sigma_{\theta_0}$ is one-to-one. $\square$

**Acknowledgments.** The authors are grateful to the anonymous referees and an Associate Editor for their valuable suggestions, which helped us greatly in improving successive versions of this paper.

J.-F. Dupuy
Laboratoire de Statistique
  et Probabilités
Université Paul Sabatier
118 route de Narbonne
31062 Toulouse cedex 9
France
E-mail: dupuy@cict.fr

I. Grama
Université de Bretagne-sud
Campus de Tohannic
56000 Vannes
France
E-mail: ion.grama@univ-ubs.fr

M. Mesbah
Laboratoire de Statistique
  Théorique et Appliquée
Université Paris 6
175 rue du Chevaleret
75013 Paris
France
E-mail: mesbah@ccr.jussieu.fr